# SEARCH COST FOR A NEARLY OPTIMAL PATH IN A BINARY TREE[1]

By Robin Pemantle

*University of Pennsylvania*

Consider a binary tree, to the vertices of which are assigned independent Bernoulli random variables with mean $p \leq 1/2$. How many of these Bernoullis one must look at in order to find a path of length $n$ from the root which maximizes, up to a factor of $1 - \varepsilon$, the sum of the Bernoullis along the path? In the case $p = 1/2$ (the critical value for nontriviality), it is shown to take $\Theta(\varepsilon^{-1} n)$ steps. In the case $p < 1/2$, the number of steps is shown to be at least $n \cdot \exp(\operatorname{const} \varepsilon^{-1/2})$. This last result matches the known upper bound from [*Algorithmica* **22** (1998) 388–412] in a certain family of subcases.

**1. Introduction.** This paper considers a problem in extreme value theory from a computational complexity viewpoint. Suppose that $\{S_{n,k} : n \geq 1, k \leq K(n)\}$ are random variables, with $K(n)$ growing perhaps quite rapidly. Let $M_n := \max_{k \leq K(n)} S_{n,k}$. A prototypical classical extreme value theorem takes the form $f_n(M_n) \to Z$, where convergence is to a constant or a distribution. When $K(n)$ grows rapidly with $n$, existence of a large value $S_{n,k}$ is not the same as efficiently being able to find such a value. There is a more compelling question from the computational viewpoint: what is the maximum value of $S_{n,k}$ that can be found by an algorithm in a reasonable time?

In this paper, we will consider the $2^n$ positions of particles in the $n$th generation of a binary branching random walk. Thus $K(n) = 2^n$ and $\{S_{n,k} : 1 \leq k \leq K(n)\}$ will be $\{S(v) : |v| = n\}$, where $|v|$ denotes the depth of a vertex $v$ and $S(v)$ is the sum of IID increments $X(w)$ over all ancestors $w$ of $v$. After reviewing known results on $M_n$, we will give upper and lower complexity bounds for finding a vertex $v$ at depth $n$ such that $S(v) \geq M_n - \varepsilon n$. It is allowed to query $X(w)$ for any $w$, and $v$ is considered "found" once we can evaluate $S(v)$, that is, once all ancestors of $v$ have been queried.

Received January 2007; revised November 2008.
[1]Supported in part by NSF Grants DMS-01-03635 and DMS-06-03821.

The problem as stated asks to maximize $S(v)$ over vertices of a fixed depth $n$. A closely related paper of Aldous [1] considers the problem of how quickly one can find a vertex $v$, at any depth, with $S(v) \geq n$. The main results herein are the lower complexity bounds proved in Theorems 3.3 and 3.4, with upper bounds included to illustrate when the lower bounds are sharp or nearly sharp. The organization of the paper is as follows. In the remainder of this section we set forth notation for branching random walks. Section 2 summarizes known limit laws for extreme values of branching random walk. A number of these results, such as Proposition 2.1, (2.7), and Propositions 2.2, 2.4 and 2.6, are used in the proofs of the lower complexity bounds. Section 3 states the main results, Section 4 proves the upper complexity bounds and other preliminary results, and Section 5 proves the lower complexity bounds.

*Notation.* The infinite rooted binary tree will be denoted $T$ and its root will be denoted 0. Write $v \in T$ when $v$ is a vertex of $T$ and $v \sim w$ when $v$ is a neighbor (parent or child) of $w$. Let $|v|$ denote the depth of $v$, that is, the distance from 0 to $v$. Write $v < w$ if $w$ is a descendant of $v$. By a "rooted path" or "branch" we mean a finite or infinite sequence $(\mathbf{0} = x_0, x_1, x_2, \ldots)$ of vertices with each $x_i$ being the parent of $x_{i+1}$. Our probability space supports random variables $\{X(v) : v \in T\}$ that are IID with common distribution that is Bernoulli with mean $p \leq 1/2$; in Proposition 3.1 below and parts of Section 2, we allow a more general common distribution but all other notation remains the same. Let $S(v) := \sum_{0 < w \leq v} X(w)$ denote the partial sums of $\{X(w)\}$ on paths from the root [in particular, $S(0) = 0$]. The maximal displacement $M_n$ is defined by

$$M_n = \max_{|v|=n} S(v).$$

The subtree from $v$ is the induced subgraph on $\{w \in T : w \geq v\}$, rooted at $v$. The subtree process $\{S(w) - S(v) : w \geq v\}$ has the same distribution as the original process $\{S(w) : w \in T\}$.

Our probability space must be big enough to support probabilistic search algorithms. We will not need to define these formally, but simply to bear in mind that there is a source of randomness independent of $\{X(v) : v \in T\}$, and that there is a filtration $\mathcal{F}_0, \mathcal{F}_1, \mathcal{F}_2, \ldots$ such that $\mathcal{F}_t$ is "everything we have looked at up to time $t$"; thus $X(v(t)) \in \mathcal{F}_t$, where $v(t)$ is the vertex we choose to inspect at time $t$, and $\{X(w) : w \neq v(1), \ldots, v(t)\}$ is independent of $\mathcal{F}_t$; without loss of generality, we assume $v(t) \in \mathcal{F}_{t-1}$, that is, any randomness needed to choose $v(t)$ is generated by time step $t-1$.

## 2. Classical extreme value results.

*Growth rate of $M_n$.* Along most infinite paths, the mean of the variables will be the mean, $p$, of their common distribution, but there will be exceptional paths where the $n$th partial sum is consistently greater than $pn$. Let $X_1, X_2, \ldots$ be IID with the same distribution as the variables $X(v)$ and let $S_n := \sum_{k=1}^n X_k$ denote the partial sums. By taking expectations, $\mathbb{P}(M_n \geq L) \leq 2^n \mathbb{P}(S_n \geq L)$. It was shown in the 1970s that for $p < 1/2$, this is asymptotically sharp (see Proposition 2.1 below). Converting this to a computation of the almost sure limiting value of $M_n/n$ requires the following large deviation computation that is by now quite standard; for details (see, e.g., [9], Section 1.9). This computation is valid for any common distribution $\mathcal{L}$ of the variables $\{X_n\}$ with exponential moments; for simplicity, since this is all we will need, assume $|X_1| \leq 1$.

Let $\mu$ denote the mean of the common distribution $\mathcal{L}$, and pick real numbers $c > \mu$ and $\lambda > 0$. Let $\phi(t) := \log \mathbb{E} e^{tX_1}$. By Markov's inequality,

$$\mathbb{P}(S_n \geq cn) \leq \frac{\mathbb{E} e^{\lambda S_n}}{e^{\lambda cn}} = \exp[n(\phi(\lambda) - c\lambda)].$$

It is easy to see that $\phi$ is convex and that when $c$ is less than the essential supremum of $X_1$, there is a unique $\lambda_*(c)$ such that this bound is minimized. Thus

$$\frac{1}{n} \log \mathbb{P}(S_n \geq cn) \leq \phi(\lambda_*(c)) - c\lambda_*(c)$$

and Chernoff's well-known theorem [8] states that this is asymptotically sharp:

$$\frac{1}{n} \log \mathbb{P}(S_n \geq cn) \to \mathsf{rate}(c) := \phi(\lambda_*(c)) - c\lambda_*(c)$$

as $n \to \infty$. The proof of this involves remarking that a certain exponential reweighting of the law $\mathcal{L}$ has mean $c$:

$$(2.1) \qquad \frac{d\mathcal{L}'}{d\mathcal{L}} = \frac{e^{\lambda_* x}}{\mathbb{E} e^{\lambda_* X_1}} \quad \Longrightarrow \quad \mathbb{E}' X_1 = c.$$

Note, for later use, that Markov's inequality extends to imply

$$(2.2) \qquad \mathbb{P}(S_n \geq cn + \beta) \leq \exp[n \cdot \mathsf{rate}(c) - \lambda_*(c) \cdot \beta].$$

The following proposition was proved in 1975 by Kingman using analytic methods, then by Biggins, using an embedded branching process (see also [10, 15] for an approach via subadditive ergodic theory).

PROPOSITION 2.1 ([16], Theorem 6, and [4], Theorem 3). *Let $c = c(\mathcal{L})$ denote the value such that $\mathsf{rate}(c) = -\log 2$. Then the maximum partial sums at each level of the binary tree satisfy*

$$\frac{M_n}{n} \to c(p) \tag{2.3}$$

*in probability as $n \to \infty$.*

In particular, when $\{X_n\}$ are Bernoulli($p$) for $0 < p < 1/2$, we have

$$\frac{1}{n} \log \mathbb{P}(S_n > qn) = H(p,q) + o(1), \tag{2.4}$$

where

$$H(p,q) := q \log \frac{p}{q} + (1-q) \log \frac{1-p}{1-q}. \tag{2.5}$$

Denoting $c := c(p) := c(\mathcal{L})$, we see that $c$ solves

$$\begin{aligned}c \log p + (1-c) \log(1-p) + c \log\left(\frac{1}{c}\right) \\ + (1-c) \log\left(\frac{1}{1-c}\right) + \log 2 = 0.\end{aligned} \tag{2.6}$$

Also, (2.2) becomes

$$\mathbb{P}(S_n \geq c(p)n + \beta) \leq 2^{-n} \exp(-\lambda_*(p)\beta). \tag{2.7}$$

*Second order behavior of $M_n$.* Fix a bounded law $\mathcal{L}$ and let $c := c(\mathcal{L})$. More accurate large deviation bounds show that $\mathbb{P}(M_n \geq cn) \to 0$, leading to two natural questions: first, estimate $\mathbb{P}(M_n \geq cn)$, and second, what correction gives the typical behavior for $M_n$? We separate into two cases, $p = 1/2$ and $p < 1/2$, which will be seen to behave rather differently in many respects.

One reason second-order results are trickier than the limit results for $M_n/n$ is that the bounds obtained by computing first moments are no longer sharp. For example, in the case of binary variables, when $p = 1/2$, the expected number of paths of length $n$ consisting entirely of ones is exactly 1. However, the actual number of such paths is the number of progeny in the $n$th generation of a critical branching process, which is known to be nonzero with probability of order $1/n$. The exact result is:

PROPOSITION 2.2 ([3], Theorem I.9.1). *For* Bernoulli(1/2) *random variables,*

$$\mathbb{P}(M_n = n) = \frac{2 + o(1)}{n}.$$

The typical behavior of $M_n$ when $p = 1/2$ will not be of great concern here; the reader may consult [6] to find a proof that $M_n = n - C \log \log n + O(1)$.

In the case $0 < p < 1/2$, the mean number of paths of length $n$ with $S(v) \geq c(p)n$ is easily shown to be $\Theta(n^{-1/2})$. The probability of existence of such a path is expected to be of order $n^{-3/2}$. Such a result has not been proved. An analogous result has, however, been proved for a branching Brownian motion. Here particles move as independent Brownian motions, each particle living for an exponential amount of time of mean one before splitting into two particles which then evolve independently. Bramson [5] shows that the maximum $M_t$ of a branching Brownian motion at time $t$ exceeds $ct$ with probability $\Theta(t^{-3/2})$, where $c = \sqrt{2}$ is the critical slope, and that $M_t = ct - \gamma \log t + O(1)$ in probability, where $\gamma = 3/(2c) = 3/2^{3/2}$. This was generalized in [7]. At slopes above the critical slope the large deviation probabilities decay exponentially: for $\lambda > \sqrt{2}$ one has $\mathbb{P}(M_t \geq \lambda t + \theta) \sim c_1(\lambda, \theta) t^{-1/2} e^{-c_2(\lambda)t}$; see [12], Theorem 6.

*Survival probability with an absorbing barrier at criticality.* For the complexity questions addressed in the present article, the crucial probabilities turn out to be *absorbing barrier* probabilities, where the events $\{M_n \geq cn\}$ and $\{M_n \geq (c-\varepsilon)n\}$ are replaced by the event that along some path from the root of length $n$ the values $S_k$ are always at least $ck$ or $(c-\varepsilon)k$, for $1 \leq k \leq n$. The term "barrier" refers to probability models in which particles are killed when they hit an absorbing barrier, which is located at $(c-\varepsilon)k$. At the critical barrier ($\varepsilon = 0$) the process dies out. Estimates of survival probabilities with a critical barrier have been published only for branching Brownian motion (though a somewhat analogous result in the discrete setting is implicit in [17], Lemma 8). Suppose each particle in a branching Brownian motion is killed when its position at any time $t$ becomes less than $ct$. Starting with a single particle at 1, Kesten estimated the tails of the survival time.

PROPOSITION 2.3 ([14], Theorem 1.3). *Consider a branching Brownian motion started with a single particle at 1, in which particles are killed when their position as a function of time becomes less than or equal to $\sqrt{2}t$ (here $\sqrt{2}$ is the critical slope). The probability for at least one particle to survive to time $t$ is $\exp(-(3\pi^2 t)^{1/3} + O(\log^2 t))$.*

REMARK. For $\lambda > \sqrt{2}$, the probabilities $\mathbb{P}(M_t \geq \lambda t)$ decay exponentially in $t$. In this regime, the quantity $\mathbb{P}(M_t \geq \lambda t)$ may be estimated up to a factor of $1 + o(1)$; such an asymptotic formula was proved in [11], Theorem 1.

*Survival probability with an absorbing barrier in the supercritical regime.* Relaxing the barrier $ck$ to the barrier $(c-\varepsilon)k$ yields a supercritical process, for which one may ask about both finite time and infinite time survival probabilities. These are the results most intimately connected with search times. The following notation is useful:

DEFINITION 1 (Survival probabilities). Let $\{X(v)\}$ be IID bounded random variables with law $\mathcal{L}$. Let $c = c(\mathcal{L})$ be the unique real number such that

$$\frac{1}{n} \log \mathbb{P}(S_n \geq cn) \to -\log 2,$$

where $S_n$ is the partial sum of $n$ IID variables with law $\mathcal{L}$. By Proposition 2.1, $M_n/n \to c$ in probability. Define the survival probability $\rho(\mathcal{L}; \varepsilon, n)$ to be the probability that there exists a path $v_0, \ldots, v_n$ of length $n$ from the root such that for all $j \leq n$, $S(v_j) \geq (c - \varepsilon)j$. In the case where $\{X(v)\}$ are Bernoulli with parameter $p$, the notation $\rho(p; \varepsilon, n)$ will be used instead of $\rho(\mathcal{L}; \varepsilon, n)$. Extend the notation to nonintegral values of $n$ by defining $\rho(\mathcal{L}; \varepsilon, n) := \rho(\mathcal{L}; \varepsilon, \lfloor n \rfloor)$.

In this notation, the quantities $\rho(p; 0, n)$ denote the tails of survival probabilities at the critical barrier. Restating Proposition 2.2, we have $\rho(1/2; 0, n) \sim 2/n$. We will be chiefly interested in the probabilities $\rho(p; \varepsilon, \infty)$ of survival to infinity once the absorbing barrier has moved so as to make the branching random walk slightly supercritical. For branching random walk with binary variables and $p = 1/2$ there is a sharp result.

PROPOSITION 2.4.

(2.8) $$\rho(\tfrac{1}{2}; \varepsilon, \infty) = \Theta(\varepsilon).$$

PROOF. Assume without loss of generality that $\varepsilon = 1/n$ for some integer, $n$. One inequality follows from the observation that a path stays above $(1 - \varepsilon)k$ for every $k < \varepsilon^{-1}$ only if it is composed entirely of ones. Therefore, from Proposition 2.2,

$$\rho(\tfrac{1}{2}; \varepsilon, \infty) \leq \rho(\tfrac{1}{2}; \varepsilon, \varepsilon^{-1} - 1) = \rho(\tfrac{1}{2}; 0, \varepsilon^{-1} - 1) \sim 2\varepsilon.$$

For the other inequality, note that $\rho(\tfrac{1}{2}; \varepsilon, \infty)$ is at least the probability that there exists an infinite path $\mathbf{0} = v_0, v_1, v_2, \ldots$, along which $X(v_i) = 1$ unless $i$ is a multiple of $n$. Let $Z_i$ count the vertices at level $i$ all of whose descendants $w$ have either $X(w) = 1$ or $n$ divides $|w|$. Then $\{Z_i\}$ are the generation sizes of a branching process that is not time-homogeneous but is periodic: the offspring generating function is $f_1(z) := (1 + z)^2/4$ at times that are not multiples of $n$ and $f_2(z) := z^2$ at times that are multiples of $n$. Using a superscript of $(k)$ to denote $k$-fold composition, we may write the generating function $\sum_k \mathbb{P}(Z_{jn} = k) z^k$ as $\Psi^{(j)}(z)$ where

$$\Psi = f_2 \circ f_1^{(n-1)}.$$

The extinction probability is the increasing limit as $j \to \infty$ of $\Psi^{(j)}(0)$. Substituting $u = 1 - z$, the survival probability is the decreasing limit of $\tilde{\Psi}^{(j)}(1)$, where $\tilde{\Psi}^{(j+1)} = g_2 \circ g_1^{(n-1)} \circ \Psi^{(j)}$ and $g_j(z) = 1 - f_j(1-z)$ for $j = 1, 2$. For $u \leq n^{-1}$ we have

$$u \geq g_1(u) = u - \frac{u^2}{4} \geq u\left(1 - \frac{1}{4n}\right)$$

and iterating $n - 1$ times gives $g_1^{(n-1)}(u) \geq (3/4)u$. Hence,

$$\tilde{\Psi}(u) \geq 2(\tfrac{3}{4}u) - (\tfrac{3}{4}u)^2 \geq u.$$

It follows that the decreasing limit of $\tilde{\Psi}^{(j)}(1)$ is at least $n^{-1}$ which is equal to $\varepsilon$, hence $\rho(\tfrac{1}{2}; \varepsilon, \infty) \geq \varepsilon$, completing the proof. $\square$

Even in the binary case, when $p < 1/2$, estimates are quite tricky. It is believed that:

CONJECTURE 1. For each $p \in (0, 1/2)$ there is a constant $\beta_p$ such that as $\varepsilon \to 0$,

$$\log \rho(p; \varepsilon, \infty) \sim -\beta_p \varepsilon^{-1/2}.$$

Furthermore, $\log \rho(p; \varepsilon, L\varepsilon^{-3/2}) \sim -\beta_{p,L} \varepsilon^{-1/2}$ with $\beta_{p,L} \to \beta_p$ as $L \to \infty$ and $\beta_{p,L} \to 0$ as $L \to 0$.

There is one subcase of the case of binary variables, for which such a result is known. Let $p_{\text{crit}}$ be the value of $p$ for which $c(p_{\text{crit}}) = 1/2$. Solving (2.6) for $p$ with $c = 1/2$ we find that

$$\tfrac{1}{2} \log p_{\text{crit}} + \tfrac{1}{2} \log(1 - p_{\text{crit}}) + \tfrac{1}{2} \log 2 + \tfrac{1}{2} \log 2 + \log 2 = 0,$$

which is equivalent to $16 p_{\text{crit}}(1 - p_{\text{crit}}) = 1$, hence $p_{\text{crit}} = (2 - \sqrt{3})/4 \approx 0.067$. Suppose we consider only pairs $(p, \varepsilon)$ such that $c(p) - \varepsilon = 1/2$. In other words, we have chosen $p$ just a little greater than $p_{\text{crit}}$ and must compute the probability that there is a path, along which, cumulatively, the ones always outnumber the zeros. Aldous showed that one may compute the probability of such an infinite path by analyzing the embedded branching process of excess ones.

PROPOSITION 2.5 ([2], Theorem 6). *For $c(p) - \varepsilon = 1/2$,*

$$\log \rho(p; \varepsilon, \infty) = -\kappa(p - p_{\text{crit}})^{-1/2} + O(1)$$

*as $p \downarrow p_0$, with*

$$\kappa = \frac{\pi \log 1/(4p_0)}{4\sqrt{1 - 2p_0}} \approx 1.11.$$

*Equivalently, since $c(p)$ has a finite derivative $\nu$ at $p_0$,*

(2.9) $$\log \rho(p; \varepsilon, \infty) = -c_* \varepsilon^{-1/2} + O(1)$$

*as $\varepsilon \to 0$ with $c(p) - \varepsilon = 1/2$, where $c_* = \kappa\sqrt{\nu}$.*

One way to prove Conjecture 1 without the restriction $c(p) - \varepsilon = 1/2$ would be to adapt Kesten's proof for Brownian motion to the random walk setting. Inspection of the nearly forty journal pages in [14] devoted to the Brownian result lead one to believe this would be possible but tedious. It is worth formulating an easier but crude result bounding $\rho(p; \varepsilon, \infty)$ from above; among other things this will clarify that the logarithms of the factors other than $\rho(p, s\varepsilon, \varepsilon^{-3/2})$ in the statement of Theorem 3.4 below are asymptotically neglible; the proof is given in Section 4.

PROPOSITION 2.6. *Fix any law $\mathcal{L}$ with mean $\mu$ supported on $[\mu - 1, \mu + 1]$. Then there is a constant, $\eta > 0$ such that for any sufficiently small $\varepsilon$,*

$$\log \rho(\mathcal{L}; \varepsilon, \infty) \leq \log \rho(\mathcal{L}; \varepsilon, \varepsilon^{-3/2}) \leq -\eta \varepsilon^{-1/2}.$$

**3. Complexity results.** An easy result, found in [13], is that a finite lookahead algorithm can produce a path with $(c(p) - \varepsilon)n$ 1's in time $g(\varepsilon)n$ for some function $g$. This suggests that we focus our effective computation question on times that are linear in $n$ and try to find the relationship between the linear discrepancy $\varepsilon$ from optimality and the linear time constant $g(\varepsilon)$.

Upper bounds for the computation time are in general easier, because finding a reasonable algorithm is easier than proving none exists. In fact, good upper bounds are obtained using a depth-first search. The notion of a depth-first search is quite standard; nevertheless, some details are required in order to avoid later ambiguities. Suppose a random set $W$ of vertices is adapted, in the sense that the event $v \in W$ is measurable with respect to $\mathcal{F}(v)$, the $\sigma$-field generated by the values $X(w)$ at all vertices $w \leq v$. A depth-first search for an infinite descending path in $W$ is the following algorithm. Label the two children of $v$ by $v0$ and $v1$, so vertices are labeled by finite sequences of zeros and ones. Order the vertices lexicographically. At time 1, examine the root; if $\mathbf{0} \notin W$ the search fails. At each subsequent time, examine the leftmost vertex $v$ (the vertex whose label is the least binary number) among children of vertices previously examined and found to be in $W$. If $v \notin W$ and is composed of all 1's, then the search fails, otherwise the search continues. Properties of the depth-first search include the following:

1. The set of examined vertices is always a subtree.
2. The sequence of examined vertices is in lexicographic order.

3. If the search continues for infinite time, then the set of vertices found to be in $W$ will contain a unique infinite path (this follows from the previous property).

Specialize now to the set $W = W_\varepsilon$ defined to be the set of vertices $v$ such that for all $w \leq v$, $S(w) \geq (c-\varepsilon)|w|$. Finding a path from the root of length $n$ in $W_\varepsilon$ is one way to locate a witness to $M_n \geq (c-\varepsilon)n$. Although there may be many witnesses outside $W_\varepsilon$, they are hard to find, so searching $W_\varepsilon$ turns out to be a pretty good way to test whether $M_n \geq (c-\varepsilon)n$. The only drawback is that the search may fail. Therefore, we define the *iterated depth-first search with parameter* $\varepsilon$, denoted by $\mathsf{IDFS}(\varepsilon)$, as follows. Recall that the subtree process is defined as the set $\{S(u) - S(v) : u \geq v\}$; thus we may define a set $W_\varepsilon(v)$, which is the set $W_\varepsilon$ of the subtree process from $v$, to be the set of $u \geq v$ such that for $v \leq z \leq u$, $S(z) - S(v) \geq (c-\varepsilon)(|z| - |v|)$.

IDFS:
Repeat until failing to terminate:

> Let $v$ be the leftmost among vertices of minimal depth that have not yet been examined, and execute a depth-first search for an infinite path in $W_\varepsilon(v)$.

Thus the algorithm begins with a depth-first search for an infinite path in $W_\varepsilon(\mathbf{0})$. If this goes on forever, then this is the whole IDFS. Otherwise, at each termination, the search begins again from a vertex none of whose descendants has been examined. Therefore, the probability of success after each termination is $\rho(\mathcal{L}; \varepsilon, \infty) > 0$. It follows that one plus the number of terminations is a geometric random variable with mean $\rho(\mathcal{L}; \varepsilon, \infty)^{-1}$ and in particular, will be finite, hence IDFS will always find an infinite path in $W_\varepsilon(v)$ for some $v$.

The next proposition uses a depth-first search to give a general upper bound in terms of certain survival probabilities; the proof is given at the beginning of the next section. The result was known to Aldous [1], though not proved in this form. For this result, binary random variables are not required.

PROPOSITION 3.1. *Let $\{X(v)\}$ be IID with any bounded distribution $\mathcal{L}$. Fix any $r < 1$ and $\varepsilon > 0$. As $n \to \infty$, the probability goes to 1 that $\mathsf{IDFS}(r\varepsilon)$ finds a vertex $v$ with $|v| = n$ and $S(v) \geq (c-\varepsilon)n$. The time it takes to do this is at most $\rho(\mathcal{L}; r\varepsilon, \infty)^{-1} n + o(n)$ in probability.*

REMARK. The appearance of $\rho$ in this bound explains why the quantities $\rho(p; \varepsilon, \infty)$ are relevant to the complexity problem.

The upper bound in the critical case follows directly from this proposition.

COROLLARY 3.2 (Upper complexity bound when $p = 1/2$). *Let $p = 1/2$. There is a $C > 0$ and an algorithm which produces a path of length $n$ having at least $(1 - \varepsilon)n$ 1's, in time at most $Cn\varepsilon^{-1}$, with probability tending to 1 as $n \to \infty$.*

PROOF. By (2.8) of Proposition 2.4, we know that $\rho(1/2; r\varepsilon, \infty) \geq c_1 r\varepsilon$ for some $c_1 \geq 0$. By Proposition 3.1, for any $\delta > 0$ there is an $n_0$ such that for $n > n_0$, IDFS($r\varepsilon$) produces the desired path by time $((c_1 r\varepsilon)^{-1} + \delta)n$ with probability at least $1 - \delta$. This proves the lemma for any $C > (c_1 r\varepsilon)^{-1}$. □

The first main result of this paper is the corresponding lower bound.

THEOREM 3.3 (Lower complexity bound when $p = 1/2$). *Let $p = 1/2$. For any search algorithm (see the discussion at the end of Section 1), for any $\kappa < 1/2$, and for all sufficiently small $\varepsilon$ (depending on $\kappa$), the probability of finding a path of length $n$ from the root with at least $(1 - \varepsilon)n$ 1's by time $\kappa\varepsilon^{-1}n$ is $O(1/n)$, uniformly in the search algorithm.*

When $p < 1/2$, lack of understanding of $\rho(p; \varepsilon, \infty)$ prevents us from stating an upper bound beyond what is inherent in Proposition 3.1. In the special case that $c(p) - \varepsilon = 1/2$, we may put Proposition 2.5 together with Proposition 3.1 to see that IDFS finds a witness to $M_n \geq (c(p) - \varepsilon)n$ by time $n \exp(C\varepsilon^{-1/2})$ for some $C > 0$. If Conjecture 1 is true, then for all $p$ and $\varepsilon$ the IDFS is likely to succeed in time $O(n \exp(C\varepsilon^{-1/2}))$. The second main result of this paper is a corresponding lower complexity bound. Because this is stated in terms of $\rho$ it is a reasonably sharp converse to Proposition 3.1.

THEOREM 3.4 (Lower complexity bound when $p < 1/2$). *Fix $p \in (0, 1/2)$ and $s > 1$. For any algorithm, the probability of finding a path of length $n$ with at least $(c(p) - \varepsilon)n$ 1's by time*

$$\frac{s-1}{4(1-c(p))} \varepsilon^{11/2} \rho(p, s\varepsilon, \varepsilon^{-3/2})^{-1} n$$

*is $O(\varepsilon^{-1} n^{-1})$.*

REMARKS. If the asymptotics for $\rho$ are as expected, then one could take $s = 1 + o(1)$ as $\varepsilon \to 0$ in such a way that

$$\log\left[\frac{s-1}{4(1-c(p))} \varepsilon^{11/2} \rho(p; s\varepsilon, \varepsilon^{-3/2})^{-1} n\right] \sim \log[\rho(p; \varepsilon, \varepsilon^{-3/2})^{-1} n].$$

This would require a regularity result on $\rho$ which is not proved. Note also, that it is expected (Conjecture 1) that

$$\log \rho(p; \varepsilon, L\varepsilon^{-3/2}) \sim -C_L \varepsilon^{-1/2},$$

but that the constant should depend on $L$, so Theorem 3.4 is at best sharp up to a constant factor in the logarithm. Finally, we note that as $\varepsilon \to 0$ with $n$ fixed, the search time ceases to grow once $\varepsilon < 1/n$, which is reflected in the fact that the probability upper bound $C\varepsilon^{-1}n^{-1}$ becomes uniformative when $\varepsilon$ is this small.

## 4. Proofs of preliminary results and upper complexity bounds.

PROOF OF PROPOSITION 3.1. Say that a vertex $v$ is *good* if there is an infinite descending path $x_0, x_1, x_2, \ldots$ from $v$ (a path where each $x_{j+1}$ is a child of $x_j$) such that $S(x_n) - S(v) \geq (c - r\varepsilon)n$ for all $n$. Such a path is called a *good path*. If IDFS($r\varepsilon$) ever examines a good vertex $v$, then it will never leave the subtree of $v$. Not every vertex on a good path is necessarily good. However, if the search algorithm encounters infinitely many good vertices $v(1), v(2), \ldots$, then, since each must be in the subtree of the previous one, these must form a chain of descendants and the sequence $v(t): t \geq 1$ must converge to a single end of the tree (an infinite descending path).

Since each vertex examined by the depth-first search has no descendants previously examined, we have

$$\mathbb{P}(v(t) \text{ is good } | \mathcal{F}_{t-1}) = \rho(\mathcal{L}; r\varepsilon, \infty)$$

for all $t$. By the conditional Borel–Cantelli lemma (e.g., [9], Theorem 4.4.11), the number of good vertices among $v(1), \ldots, v(n)$ is almost surely $\rho(r\varepsilon)n + o(n)$. Hence, after the time $\tau_n$ that the $n$th good vertex is examined, the path from $v(\tau_1)$ to $v(\tau_n)$ has the property that any vertex $w$ on the path has

$$S(w) - S(v(\tau_1)) \geq (c - r\varepsilon)(|w| - |v(\tau_1)|).$$

Recalling that $r < 1$, we see there is a random $N$ such that for all vertices $v$ on the infinite path chosen by the algorithm, if $|v| \geq N$ then $S(v) \geq (c - \varepsilon)|v|$. The conclusion of the proposition follows. □

By Brownian scaling, for a mean-zero, finite variance random walk $\{S_n\}$, we have $\log \mathbb{P}(S_1, \ldots, S_n \in [-L, L]) \sim -Cn/L^2$. It is convenient to record a lemma giving an explicit constant for the upper bound, uniform over all walks with a given variance.

LEMMA 4.1. *Let $\{S_n\}$ be a random walk whose increments are bounded by 1 and have mean zero and variance $\sigma^2 > 0$. Then for $L \geq 1$, the probability of the walk staying in an interval $[-L, L]$ up to time $N$ is bounded above by*

$$\mathbb{P}(S_1, \ldots, S_N \in [-L, L]) \leq \exp\left(-\frac{\sigma^2 N}{36eL^2}\right),$$

*provided that the exponent is less than $-1/4$, that is, $N > 9eL^2/\sigma^2$.*

PROOF. For any $n \leq N$, the event that $S_1, \ldots, S_N \in [-L, L]$ implies that for each $j \leq k \leq j+n$, $|S_k - S_j| \leq 2L$. Breaking into $\lfloor N/n \rfloor$ time blocks of size $n$, plus a possible leftover segment, independence of the increments implies that

$$(4.1) \qquad \mathbb{P}(S_1, \ldots, S_N \in [-L, L]) \leq \mathbb{P}(S_1, \ldots, S_n \in [-2L, 2L])^{\lfloor N/n \rfloor}.$$

Later, we will choose

$$(4.2) \qquad n = \left\lceil \frac{8eL^2}{\sigma^2} \right\rceil.$$

For now, we let $n$ and $\alpha$ be arbitrary and we let $\tau_\alpha := \inf\{k : |S_k| \geq \alpha\sqrt{n}\}$ be the time for the random walk to exit the interval $[-\alpha\sqrt{n}, \alpha\sqrt{n}]$. Let us obtain an upper bound on $\mathbb{P}(\tau_\alpha > n)$. Clearly, $\mathbb{E}S^2_{\tau_\alpha \wedge (n+1)} \leq (\alpha\sqrt{n}+1)^2$ because $S_{\tau_\alpha \wedge (n+1)} \in [-1-\alpha\sqrt{n}, 1+\alpha\sqrt{n}]$. Hence,

$$(\alpha\sqrt{n}+1)^2 \geq \mathbb{E}S^2_{\tau_\alpha \wedge (n+1)} \geq \sum_{j=0}^{n} \sigma^2 \mathbb{P}(\tau_\alpha > j) \geq \sigma^2 (n+1) \mathbb{P}(\tau_\alpha > n).$$

Choosing $\alpha = \sigma/\sqrt{2e}$ and using $(a+c)^2 \leq 2(a^2 + c^2)$ now gives

$$(4.3) \quad \mathbb{P}(\tau_{(2e)^{-1/2}\sigma} > n) \leq \frac{((2e)^{-1/2}\sigma\sqrt{n}+1)^2}{\sigma^2 n} \leq e^{-1} + \frac{2}{\sigma^2 n} \leq e^{-1/2},$$

once $\sigma^2 n \geq 8e$. Now, choosing $n$ as in (4.2) implies that $\alpha\sqrt{n} = \sigma\sqrt{n}/\sqrt{2e} \geq 2L$ and hence by (4.1) and (4.3),

$$\mathbb{P}(S_1, \ldots, S_N \in [-L, L]) \leq \exp\left(-\frac{1}{2}\left\lfloor \frac{N}{n} \right\rfloor\right).$$

The proof is finished by observing that $\lceil n \rceil \leq 9eL^2/\sigma^2$ (because $L^2 > 1 > \sigma^2$) and hence that

$$\left\lfloor \frac{N}{n} \right\rfloor \geq \left\lfloor \frac{N\sigma^2}{9eL^2} \right\rfloor \geq \frac{N\sigma^2}{18eL^2},$$

once $N\sigma^2 \geq 9eL^2$. □

PROOF OF PROPOSITION 2.6. We need to find a constant $\eta > 0$ such that

$$(4.4) \qquad \rho(\mathcal{L}; \varepsilon, \varepsilon^{-3/2}) < \exp(-\eta \varepsilon^{-1/2}).$$

This is a standard "squeezing" argument: the probability space is broken into two parts. One has small measure because some particle is found at a position that is greater by $\alpha \varepsilon^{-1/2}$ than it should be, for some constant $\alpha > 0$; conditioning on the complement of this event squeezes the path below $c(\mathcal{L})k + \alpha \varepsilon^{-1/2}$, but above $(c(\mathcal{L}) - \varepsilon)k$, at each level $k$; the chance of a random walk trajectory remaining in such a tube is small enough to make the expected number of such trajectories small.

To make this precise, begin by letting $c = c(\mathcal{L})$ and $\lambda_*$ be as in Proposition 2.1. For any positive integer $N$, denote by $G(N, \varepsilon, \alpha)$ the event

$$\{\exists v : |v| \leq N \text{ and } S(v) \geq c(\mathcal{L})|v| + \alpha\varepsilon^{-1/2}\}.$$

Applying (2.2) to $S(v)$ for each of $2^n$ vertices $v$ at each generation $n \leq N$, we see that

(4.5) $$\mathbb{P}[G(N, \varepsilon, \alpha)] \leq N \exp(-\lambda_* \alpha \varepsilon^{-1/2}).$$

Next, set $N := \alpha \varepsilon^{-3/2}$ where $\alpha \leq 1$ is a positive parameter that will be specified later. Let $\mathcal{L}$ denote the common law of the Bernoulli($p$) variables $\{X_n\}$, let $\mathcal{L}'$ denote the law of a Bernoulli($c(\mathcal{L})$) variable and $\mathcal{L}''$ denote the law of a compensated Bernoulli($c(\mathcal{L})$) variable. Let $\mathbf{Q}$ (resp., $\mathbf{Q}', \mathbf{Q}''$) denote the law of a sequence $\{X_n\}$, that is, IID with law $\mathcal{L}$ (resp., $\mathcal{L}', \mathcal{L}''$). Let $\sigma^2$ be the common variance of $\mathcal{L}'$ and $\mathcal{L}''$ and let $c_1$ denote the constant $-\sigma^2/(36e)$ from Lemma 4.1. Let $H$ denote the event that $|S_n - c(\mathcal{L})n| \leq \alpha\varepsilon^{-1/2}$ for $1 \leq n \leq N$. Applying Lemma 4.1 to the law $\mathcal{L}''$ we see that

(4.6) $$\mathbf{Q}'(H) = \mathbf{Q}''(|S_n| \leq \alpha\varepsilon^{-1/2} \text{ for all } n \leq N)$$
$$\leq \exp\left(-\frac{c_1}{\alpha}\varepsilon^{-1/2}\right).$$

For any measures $\nu$ and $\pi$ and any event $A$, we have $\nu(A) \leq \pi(A) \cdot \sup_{\omega \in A} (d\nu/d\pi)(\omega)$. We may therefore use (2.1) to convert (4.6) into an estimate for $\mathbf{Q}(H)$: plugging $\beta = -\alpha\varepsilon^{-1/2}$ and $\mathsf{rate}(c) = \log(1/2)$ into (2.1) gives

$$\mathbf{Q}(H) \leq \sup_{x \geq c(\mathcal{L})n - \alpha\varepsilon^{-1/2}} \frac{\mathbb{E}e^{\lambda^* X_1}}{e^{\lambda_* x}} \mathbf{Q}'(H)$$
$$\leq 2^{-N} \exp(\lambda_* \alpha \varepsilon^{-1/2}) \exp\left(-\frac{c_1}{\alpha}\varepsilon^{-1/2}\right).$$

As $\alpha \downarrow 0$, the quantity $\lambda_* \alpha - c_1/\alpha$ converges to $-\infty$, hence we may pick $\alpha \in (0, 1)$ such that $\lambda_* \alpha - c_1/\alpha \leq -\lambda_* \alpha$. Fixing this value of $\alpha$ and denoting $\eta := \lambda_* \alpha$, we have

$$\mathbf{Q}(H) \leq 2^{-N} \exp(-\eta \varepsilon^{-1/2}).$$

We apply this to the variables $\{S(w) : w \leq v\}$ for the branching random walk on the binary tree, where $v$ is any vertex at depth $N$. There are $2^N$ such vertices, whence the probability that some path $\mathbf{0} = v_0, \ldots, v_N$ satisfies $|S(v_n) - c(\mathcal{L})n| \leq \alpha\varepsilon^{-1/2}$ for all $n \leq N$ is bounded above by $\exp(-\eta\varepsilon^{-1/2})$. Combining this with (4.5) shows that

$$\rho(\mathcal{L}; \varepsilon, \alpha\varepsilon^{-3/2}) \leq (N+1)\exp(-\eta\varepsilon^{-1/2}).$$

Choosing a slightly smaller value of $\eta$, we may absorb the factor of $N + 1$. Because $\alpha$ is at most 1 and $\rho(\mathcal{L}; \varepsilon, N)$ is decreasing in $N$, the proof is complete. $\square$

**5. Proofs of lower complexity bounds.** An easy lemma needed at the end of each of the two proofs is the following:

LEMMA 5.1.  *Let $\{X_t : t = 1, 2, 3, \ldots\}$ be adapted to the filtration $\{\mathcal{F}_t\}$ and have partial sums $S_t := \sum_{k=1}^{t} X_k$. Suppose there are numbers $\beta_t$ and $\alpha_t$ such that for all $t$,*

$$\mathbb{E}(X_t \mid \mathcal{F}_{t-1}) \leq \beta_t;$$
$$\mathbb{E}(X_t^2 \mid \mathcal{F}_{t-1}) \leq \alpha_t.$$

*Let $\overline{\beta_t} := t^{-1} \sum_{s=1}^{t} \beta_s$ and $\overline{\alpha_t} := t^{-1} \sum_{s=1}^{t} \alpha_s$. Then for any $T$ and any $\beta' > \overline{\beta_T}$,*

$$\mathbb{P}(S_T > T \beta') \leq \frac{\overline{\alpha_T}}{(\beta' - \overline{\beta_T})^2} T^{-1}.$$

*In the special case $\beta_t \equiv \beta, \alpha_t \equiv \alpha$, this becomes*

(5.1) $$\mathbb{P}(S_T > T \beta') \leq \frac{\alpha}{(\beta' - \beta)^2} T^{-1}.$$

PROOF.  Let $\mu_t = \sum_{k=1}^{t} \mathbb{E}(X_k \mid \mathcal{F}_{k-1})$. Then $\{S_t - \mu_t\}$ is a martingale and

$$\mathbb{E}(S_t - \mu_t)^2 = \sum_{k=1}^{t} \mathrm{Var}(X_k \mid \mathcal{F}_{k-1}) \leq \sum_{k=1}^{t} \mathbb{E}(X_k^2 \mid \mathcal{F}_{k-1}) \leq t \overline{\alpha_t}.$$

Using the inequality $\mu_t \leq t \overline{\beta_t}$, we then have, by Chebychev's inequality,

$$\mathbb{P}(S_T > \beta' T) \leq \mathbb{P}(S_T - \mu_T \geq (\beta' - \overline{\beta_T}) T) \leq \frac{T \overline{\alpha_T}}{(\beta' - \overline{\beta_T})^2 T^2}. \quad \square$$

PROOF OF THEOREM 3.3.  It suffices to prove the result when $\varepsilon = 1/(2b)$ is the reciprocal of an even integer and $n$ is even. For such values of $\varepsilon$ and $n$, divide the vertices of $T$ into two classes. Label $v$ as *good* if there is a path of length $b$ descending from $v$ on which the labels are all equal to one; label all other vertices *bad*. Suppose $\gamma = (x_0, \ldots, x_{n-1})$ is a path of length $n - 1$ from the root and that at most $\varepsilon n$ vertices $v \in \gamma$ have $X(v) = 0$. Then at most $b \varepsilon n + b$ vertices $v \in \gamma$ are bad, because if $x_j$ is bad for $j < n - b$ then at least one of $x_j, \ldots, x_{j+b}$ must be labeled with a zero. It follows that least $n/2 - b$ of the vertices in $\gamma$ are good.

Say that our search algorithm *does not jump* if each successive vertex inspected is a neighbor of the root or of a vertex previously inspected. These algorithms have the property that whenever you peek at a vertex you know nothing about its descendant tree.

CONJECTURE 2. No algorithm finds a path with at least $(1-\varepsilon)n$ 1's in a shorter average time than the best algorithm that does not jump. (See [1], Conjecture 5.1, for a similar conjecture.)

If the conjecture is true, then the proof of the theorem is very short: each new vertex we peek at has probability $O(1/n)$ of being good, independent of the past; in time $o(n^2)$, we can therefore find only $o(n)$ good vertices, and in particular, we cannot find $n/2$ good vertices. In absence of a proof of the conjecture, the proof of the theorem continues as follows.

Given a search algorithm producing a sequence $\{v(t) : t \geq 1\}$ of examined vertices, define sets $A(t)$ as follows. A vertex $x$ is in $A(t)$ if all of the following hold:

(i) $x \notin \bigcup_{s<t} A(s)$;
(ii) $x = v(t)$ or $x$ is an ancestor of $v(t)$ and $|x| > |v(t)| - b$;
(iii) there is a descending path from $x$ of length $b$, passing through $v(t)$, all of whose vertices $w$ have $X(w) = 1$.

In other words, $x \in A(t)$ if $t$ is the first time a vertex $v$ is peeked at that lies on a path of length $b$ of 1's descending from $x$. Think of $A(t)$ as an accounting scheme which marks good vertices as "found" as soon as their subtree is explored. To avoid confusion, note that $A(t)$ is not measurable with respect to $\mathcal{F}_t$: good vertices "found" at time $t$ are not known to be good until much later. If a path $\gamma = (x_0, \ldots, x_{n-1})$ has at most $\varepsilon n$ zeros on it, and this whole path has been found by our search algorithm by time $t$, then there are at least $n/2 - b$ values of $j$ such that $X(x_j) = \cdots = X(x_{j+b-1}) = 1$. For these values of $j$, the vertex $x_j$ is good and is in $A(s)$ for some $s \leq t$. Thus finding $\gamma$ by time $t$ implies

$$(5.2) \qquad \left| \bigcup_{s \leq t} A(s) \right| \geq \frac{n}{2} - b.$$

Now we bound the conditional mean and variance of $|A(t)|$ given $\mathcal{F}_{t-1}$. Let $y_j(t)$ denote the ancestor of $v(t)$ going back $j$ generations. The possible elements of $A(t)$ are $v(t) = y_0(t), y_1(t), \ldots, y_{b-1}(t)$. The event $y_j(t) \in A(t)$ is contained in the intersection of the events $G_j := \{X(y_j(t)) = \cdots = X(y_0(t)) = 1\}$, $G'_j := \{y_i \notin A(s) \ \forall 0 \leq i \leq j, 0 < s < t\}$ and the event $H_j$ that if $j < b-1$ then there is a descending path of length $b-1-j$ from a child of $v(t)$ labeled by 1's and disjoint from $\bigcup_{s<t} A(s)$. Clearly $G'_j \in \mathcal{F}_t$ and on $G'_j$, $\mathbb{P}(G_j \mid \mathcal{F}_{t-1}) = 2^{-j-1}$. Also, on $G'_j$, Proposition 2.2 and the definition of $\rho(p, \varepsilon, n)$ implies that for $j < b - 1$,

$$\mathbb{P}(H_j \mid \mathcal{F}_{t-1}, G_j) = \rho(1/2, 1, b-1-j) \sim \frac{2}{b-j}$$

as $b - j \to \infty$, independently of $t$. Putting this together gives

$$\mathbb{E}(|A(t)| \mid \mathcal{F}_{t-1}) \leq \sum_{j=0}^{b-1} 2^{-j-1} \frac{2 + u(b-1-j)}{b-j}$$

for some function $u$ tending to zero. Since $2/b = 4\varepsilon$, the $j$th term is $2^{-j-1}(4\varepsilon + o(\varepsilon))$ uniformly in $t$ and summing the last expression gives

$$\mathbb{E}(|A(t)| \mid \mathcal{F}_{t-1}) \leq 4\varepsilon + o(\varepsilon),$$

uniformly in $t$ as $\varepsilon \to 0$.

To bound the second moment, compute

$$\mathbb{E}(|A(t)|^2 \mid \mathcal{F}_{t-1}) \leq \sum_{0 \leq j,k \leq b-1} \mathbb{P}(y_j, y_k \in A(t) \mid \mathcal{F}_{t-1})$$

$$\leq \sum_{j=0}^{b-1} (2j+1) \mathbb{P}(y_j \in A(t) \mid \mathcal{F}_{t-1}).$$

The probability on the right-hand side is the probability each of the vertices $y_0, \ldots, y_j$ being marked with a 1, and simultaneously, of the existence of a path of length $b - 1 - j$ of vertices descending from $v(t)$ also all bearing 1's. This probability is equal to

$$\sum_{j=0}^{b-1} \frac{2j+1}{2^{j+1}} \frac{2 + u(b-1-j)}{b-j},$$

which is asymptotic to $4\varepsilon \sum_{j=0}^{\infty} (2j+1)/2^{j+1} = 12\varepsilon$ as $\varepsilon \to 0$. Now fix $\kappa < 1/2$ and use Lemma 5.1 with $X_t = |A(t)|$, $\beta = 4\varepsilon + o(\varepsilon)$, $\alpha = 12\varepsilon + o(\varepsilon)$, $T = \kappa n \varepsilon^{-1}/4$ and

$$\beta' = \frac{n/2 - b}{T} = 4\frac{n/2 - \varepsilon^{-1}/2}{\kappa n \varepsilon^{-1}} = 4\frac{1/2 - \varepsilon^{-1}/(2n)}{\kappa \varepsilon^{-1}} = \frac{4\varepsilon}{2\kappa} + O\left(\frac{1}{n}\right),$$

uniformly in $\varepsilon$. The conclusion, recalling (5.2), is that the probability of finding a path $\gamma$ of length $n$ with at most $\varepsilon n$ zeros on it by time $T$, is at most

$$\mathbb{P}\left(\left|\bigcup_{s \leq 4\kappa n \varepsilon^{-1}} A(s)\right| \geq \frac{n}{2} - b\right) \leq \mathbb{P}\left(\sum_{s \leq 4\kappa n \varepsilon^{-1}} |A(s)| \geq \frac{n}{2} - b\right)$$

(5.3)
$$\leq \theta(\varepsilon, n) n^{-1},$$

where
$$\theta(\varepsilon, n) := \frac{\alpha}{(\beta' - \beta)^2}(nT^{-1}).$$

Computing $\theta(\varepsilon, n)$ we have
$$\beta' - \beta = 4\varepsilon\left(\frac{1}{2\kappa} - 1\right) + o(\varepsilon) + O\left(\frac{1}{n}\right)$$

and hence
$$\theta(\varepsilon, n) = 4\frac{(12 + o_\varepsilon(1))\varepsilon}{\varepsilon^2(1/(2\kappa) - 1 + o_\varepsilon(1) + O(1/(n\varepsilon)))^2}\kappa^{-1}\varepsilon$$
$$\to \frac{48\kappa}{(1/2 - \kappa)^2}$$

as $\varepsilon \to 0$ and $n\varepsilon \to \infty$. In particular, for sufficiently small $\varepsilon > 0$ and all $n$, we see that $\theta(\varepsilon, n)$ is bounded and the conclusion of the theorem follows from (5.3). □

PROOF OF THEOREM 3.4. It suffices to prove the theorem when $\varepsilon = b^{-2/3}$ for some integer $b$. Fix $p < 1/2$ and $s > 1$. The strategy is again to show that one must find a lot of good vertices and that a good vertex is hard to find. This time, the set of good vertices is the set $R(p; \varepsilon, b)$ of vertices $v$ for which there is a descending path $x_0, \ldots, x_b$ from $v$ such that for each $1 \leq j \leq b$,
$$\sum_{i=1}^{j} X(x_i) \geq (c(p) - \varepsilon) \cdot j.$$

Observe that $\mathbb{P}(v \in R(p; \varepsilon, b)) = \rho(p; \varepsilon, b)$ for all $v$.

LEMMA 5.2 (Must find good vertices). *Let $s > 1$ and suppose that $\gamma = \gamma(\varepsilon, n) = (x_0, \ldots, x_n)$ is a path of length $n$ from the root with at least $(c(p) - \varepsilon) \cdot n$ ones. Then there are $0 < \varepsilon_0 < n_0 < \infty$ such that for $\varepsilon \leq \varepsilon_0$ and $n \geq n_0$, the number of vertices in $\gamma(\varepsilon, n) \cap R(p, s\varepsilon, b)$ is at least*
$$\left\lfloor \frac{s-1}{2(1-c(p))} n\varepsilon^{5/2} \right\rfloor.$$

REMARK. Again, the proof finds this many good vertices that are not only elements of $\gamma \cap R$, but for which the values of $X(v)$ for $v \in \gamma$ are a witness to this.

PROOF OF LEMMA 5.2. Color the vertices of $\gamma$ red and blue under the following recursive rule. Let $\tau_0 = 0$. Recursively define $\tau_{j+1}$ to be $\tau_j + k$ where $k$ is the least positive integer less than $b$ for which

$$S(x_{\tau_j+k}) \leq S(x_{\tau_j}) + k(c(p) - s\varepsilon),$$

if such an integer exists, and is equal to $b$ otherwise. Let $J$ be the least $j$ for which $\tau_j \geq n$. All vertices in the list $x_{\tau_j+1}, \ldots, x_{\tau_{j+1}}$ receive the same color. The color is red if $\tau_{j+1} < n$ and $\tau_{j+1} < \tau_j + b$ and blue otherwise. Denote the set of red and blue vertices by red and blue, respectively; see Figure 1 for an example of this.

The sum $S_{\mathsf{red}} := \sum_{v \in \mathsf{red}} X(v)$ is equal to the sum over all $j$ for which $\tau_{j+1} < n \wedge (\tau_j + b)$ of $S(\tau_{j+1}) - S(\tau_j)$. The sum over each such segment of $X(v)$ is at most $(\tau_{j+1} - \tau_j)(c(p) - s\varepsilon)$, when

$$S_{\mathsf{red}} \leq |\mathsf{red}|(c(p) - s\varepsilon).$$

On the other hand, for each $j$ such that the vertices $x_{\tau_j+1}, \ldots, x_{\tau_{j+1}}$ are blue, either $\tau_j > n - b$ or $x_{\tau_j} \in \gamma \cap R(p; s\varepsilon, b)$. Thus the number of blue vertices is at most $b(|\gamma \cap R(p; s\varepsilon, b)| + 1)$. Using $|\mathsf{red}| + |\mathsf{blue}| = n$ and $S_{\mathsf{blue}} \leq |\mathsf{blue}|$, we have the inequalities

$$(c(p) - \varepsilon)n \leq S(x_n)$$
$$\leq (n - |\mathsf{blue}|)(c - s\varepsilon) + |\mathsf{blue}|$$
$$= n(c(p) - s\varepsilon) + |\mathsf{blue}|(1 - c(p) + s\varepsilon)$$
$$\leq n(c(p) - s\varepsilon) + (b(|\gamma \cap R(p; s\varepsilon, b)| + 1))(1 - c(p) + s\varepsilon).$$

Solving and plugging in $b = \varepsilon^{-3/2}$ yields

$$|\gamma \cap R(p; s\varepsilon, b)| \geq n\varepsilon^{5/2} \left[ \frac{s-1}{1 - c(p) + s\varepsilon} - \frac{1}{n\varepsilon} \right].$$

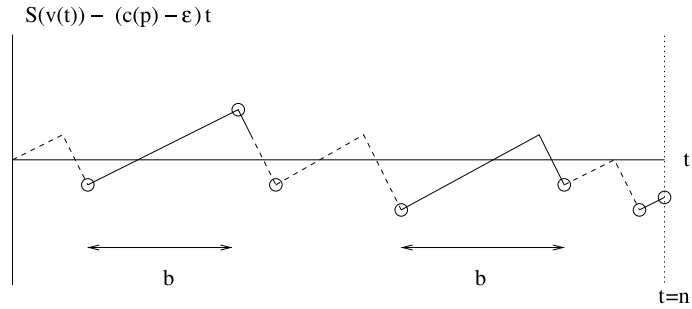

FIG. 1.   *Times $\tau_j$ are marked by hollow dots; red segments are dashed, blue segments are solid.*

When $\varepsilon$ is sufficiently small and $n\varepsilon$ is sufficiently large, the quantity in square brackets is at least half of $(s-1)/(1-c(p))$, as desired. This finishes the proof of the lemma because the result is trivial when $\varepsilon^{5/2} n < 2(1-c(p))/(s-1)$. □

CONTINUATION OF PROOF OF THEOREM 3.4. Because we do not care about factors that are polynomial in $\varepsilon$, the count is not as delicate as in the proof of Theorem 3.3. Define $A(t)$ to be the set of vertices $x$ such that there is a descending path from $x$ of length $b$, passing through $v(t)$, such that for all $j \le b$, the initial segment of length $j$ has at least $(c(p) - s\varepsilon)j$ 1's, and such that $t$ is minimal for this to hold. Formally, $x \in A(t)$ if:

(i) $x \notin \bigcup_{s<t} A(s)$;
(ii) there is a descending path $x = y_0, y_1, \ldots, y_{b-1}$ from $x$ containing $v(t)$;
(iii) for all $j \le b$, $\sum_{i=0}^{j-1} X(y_i) \ge (c(p) - s\varepsilon)j$.

Again, given $\mathcal{F}_{t-1}$, the possible elements of $A(t)$ are the ancestors of $v(t)$ back $b-1$ generations. For each ancestor $y$, $\mathbb{P}(y \in A(t) \mid \mathcal{F}_{t-1})$ is bounded above by $\rho(p; s\varepsilon, b)$. Therefore,

$$\mathbb{E}(|A(t)| \mid \mathcal{F}_{t-1}) \le b\rho(p; s\varepsilon, b).$$

For the second moment, it suffices to note the upper bound:

$$\mathbb{E}(|A(t)|^2 \mid \mathcal{F}_{t-1}) \le b\mathbb{E}(|A(t)| \mid \mathcal{F}_{t-1})$$
$$\le b^2 \rho(p; s\varepsilon, b).$$

Let $N = \lfloor \frac{(s-1)}{2(1-c(p))} n\varepsilon^{5/2} \rfloor$. By Lemma 5.2, for any $T > 0$,

(5.4)
$$\mathbb{P}(\text{finding a witness to } M_n \ge (c(p) - s\varepsilon)n \text{ by time } T)$$
$$\le \mathbb{P}\left(\left|\bigcup_{s \le T} A(s)\right| \ge N\right).$$

Let $\alpha = b\rho(p; s\varepsilon, b)$, $\beta = b^2 \rho(p; s\varepsilon, b)$, $\beta' = 2\beta$ and $T = N/\beta'$. Applying (5.1) of Lemma 5.1 bounds the right-hand side of (5.4) from above by

$$\frac{\alpha}{(\beta' - \beta)^2} \frac{1}{T} = \frac{b\rho(p; s\varepsilon, b)}{b^4 \rho(p; s\varepsilon, b)^2} \frac{2b^2 \rho(p; s\varepsilon, b)}{N}$$
$$= \frac{2}{bN}.$$

This goes to zero as $n \to \infty$; in fact, $b^{-1}N^{-1} = \varepsilon^{3/2} N^{-1} = O(\varepsilon^{-1} n^{-1})$. It follows that the probability of finding a witness to $M_n \ge (c(p) - \varepsilon)n$ by time $T$ is $O(\varepsilon^{-1}n^{-1})$. Using the fact that $\lfloor x \rfloor \geq x/2$ once $x \geq 1$ the proof is completed by observing that, once $N \geq 1$,

$$\begin{aligned}T &= \frac{N}{2\beta} \geq \frac{s-1}{4(1-c(p))} \frac{n\varepsilon^{5/2}}{2\varepsilon^{-3}\rho(p;s\varepsilon,\varepsilon^{-3/2})} \\ &= \frac{s-1}{4(1-c(p))} \varepsilon^{11/2} \rho(p;s\varepsilon,\varepsilon^{-3/2})^{-1} n. \end{aligned} \qquad \Box$$

**Acknowledgments.** The questions in this paper were brought to my attention by Yuval Peres at a meeting of the I.E.S. in Banff, 2003. A great debt is owed as well to the two referees for the improvement of expositional clarity as well as the removal of errors.

DEPARTMENT OF MATHEMATICS
UNIVERSITY OF PENNSYLVANIA
209 SOUTH 33RD STREET
PHILADELPHIA, PHILADELPHIA 19104
USA
E-MAIL: pemantle@math.upenn.edu